\documentclass[12pt,a4paper]{scrartcl}

\usepackage{amsthm,amsfonts,amsmath,latexsym,mathrsfs}
\usepackage{mathtools}

\usepackage[utf8]{inputenc}

\usepackage{libertine}
\usepackage[T1]{fontenc}
\usepackage[libertine]{newtxmath}

\usepackage[pdftex,bookmarks=true]{hyperref}
\usepackage{enumitem}

\usepackage{authblk}

\setlist{noitemsep}

\theoremstyle{definition}
\newtheorem{definition}{Definition}

\theoremstyle{remark}
\newtheorem{remark}[definition]{Remark}

\theoremstyle{plain}
\newtheorem{theorem}[definition]{Theorem}
\newtheorem{result}[definition]{Result}
\newtheorem{lemma}[definition]{Lemma}
\newtheorem{proposition}[definition]{Proposition}
\newtheorem{corollary}[definition]{Corollary}

\newcommand{\eps}{\varepsilon}

\newcommand{\lam}{\lambda}

\newcommand{\symq}{\pi_{\mathrm{Sym}}}
\newcommand{\symf}[1]{#1_{\mathrm{Sym}}}

\newcommand{\OM}{\Omega}
\newcommand{\disk}{\mathbb{D}}

\newcommand{\hol}{\mathcal{O}}

\newcommand{\bcdot}{\boldsymbol{\cdot}}

\newcommand{\C}{\mathbb{C}}

\newcommand{\Z}{\mathbb{Z}}
\newcommand{\D}{\mathbb{D}}

\newcommand{\h}{\mathfrak{h}}

\newcommand{\sympoly}{\mathbb{G}^n}
\newcommand{\Sym}{\mathrm{Sym}}

\makeatletter
  \def\blfootnote{\gdef\@thefnmark{}\@footnotetext}
\makeatother

\begin{document}

\title{Schwarz lemmas via the pluricomplex Green's function}
\author{Jaikrishnan Janardhanan \thanks{Jaikrishnan Janardhanan is supported by a DST-INSPIRE fellowship from the 
Department of Science and Technology, India.} \\ 
\href{mailto:jaikrishnan@iitm.ac.in}{\normalsize jaikrishnan@iitm.ac.in} }
\affil{Department of Mathematics, Indian Institute of Technology Madras, Chennai
600036, India}
\date{\vspace{-5ex}}

\maketitle

\begin{abstract}
  We prove a version of the Schwarz lemma for holomorphic mappings from the unit
  disk into the symmetric product of a Riemann surface.
  Our proof is function-theoretic and  self-contained.  
  The main novelty in our proof is the use of the pluricomplex Green's function.
  We also prove several other Schwarz lemmas using this function.
\end{abstract}

\blfootnote{\textup{2010} \textit{Mathematics Subject Classification}:
Primary: 32H15, 32H20, 32U35}


\section{Introduction}
The main result of this article is the following: 
\begin{theorem}\label{T:main}
  Let $X$ be a Riemann surface and 
  $f:\disk \to \Sym^n(X)$ be holomorphic. Then 
  \begin{equation}\label{E:main}
    \mathcal{H}^n_{\mathcal{M}_X}\big(\symq^{-1}(f(x)), \symq^
    {-1}(f(x_0)\big) \leq \mathcal{M}_\disk(x,x_0), \ \ \forall x,x_0 \in \disk.
  \end{equation}
\end{theorem}
\noindent The notation used in the theorem will be explained in Section~\ref{S:tools}.
Briefly, $\Sym^n(X)$ is the $n$-fold symmetric product of $X$ and
$\symq$ is the natural map from $X^n$ to $\Sym^n(X)$. The
Möbius pseudodistance associated to $X$ is denoted $\mathcal{M}_X$ and 
$\mathcal{H}_{\mathcal{M}_X}$ is the Hausdorff pseudodistance induced on
subsets of $X$ by $\mathcal{M}_X$. Given a point
$p \in \Sym^n(X), \symq^{-1}(p)$ here denotes the subset of $X$ comprising the
coordinate points of (any element of) $\symq^{-1}\{p\}$. Note that in our
notation, $\symq^{-1}\{p\}$ is \emph{not} the same as $\symq^{-1}(p)$.

\begin{remark}
  Note that Theorem~\ref{T:main} is trivially true whenever $X$ is a compact
  Riemann surface. Also for a domain    
  $D \subset \C$, it is easy to see that either $D$ is Carathéodory hyperbolic 
  (i.e, the Möbius pseudodistance is a distance) or $\mathcal{M}_D \equiv 0$.
  This is \emph{not} true for Riemann surfaces; see \cite{stanton1975bounded}.
  We emphasize that Theorem~\ref{T:main} applies to \emph{all} Riemann surfaces
  including those for which the Möbius pseudodistance is \emph{not} a distance
  but yet \emph{not} identically $0$.
\end{remark}

The genesis of Theorem~\ref{T:main} is a result by Nokrane and Ransford \cite[Theorem~1.1]
{ransford2001schwarz} which is in the setting of
algebroid multifunction taking values in the unit disk. This was later extended
to proper holomorphic correspondences from the unit disk to any bounded planar 
domain by Chandel \cite[Theorem~1.7]{chandel2017spectral}. In our notation, the
result of Nokrane and Ransford is Theorem~\ref{T:main} with $X = \disk$ while
that of Chandel is the case when $X$ is any bounded planar domain. 

Our motivation for formulating and proving Theorem~\ref{T:main} comes from an
earlier work  \cite{haridas2018note} in which we investigated the Minkowski
function of a quasi-balanced domain. During the
course of our study, we realized that a special case of 
\cite[Theorem~1.1]{ransford2001schwarz} follows easily from simple
observations about the Minkowski function and an extremal function (now
popularly known as
the \emph{pluricomplex Green's function}) studied by Lempert 
\cite{lempert1981metric}, Klimek \cite{klimek1985extremal} and
Demailly \cite{demailly1987measures}. A natural
question to ask is whether these elementary observations have wider
applicability.

The symmetrized bidisk and polydisk have been the subject of intense research
for the past two decades; see, for instance,
\cite{agler2001schwarz,edigarian2005geometry,nikolovsymmetrized,agler2013extremal}. 
More recently, the symmetric product of more general objects has also been
studied by several researchers 
\cite{chakrabarti2015function,bharali2018proper,chakrabarti2018proper,zwonek2018function}. The symmetric product of a Riemann surface can be given a
natural complex structure that makes it into a complex manifold. It is, therefore, natural to look for an extension of
the original result of Nokrane and Ransford in the setting of symmetric
products of a Riemann surface and
Theorem~\ref{T:main} is the desired extension. The proofs of both \cite
[Theorem~1.1]{ransford2001schwarz} and
\cite[Theorem~1.7]{chandel2017spectral} rely on the holomorphic functional 
calculus (the underlying 
Banach algebra being the space of $n \times n$ complex matrices). It is unclear
how these ideas can be generalized to the setting of an arbitrary Riemann
surface.
In contrast, our proof of Theorem~\ref{T:main} is almost entirely 
self-contained and uses tools solely from complex analysis. Specifically, we 
require only basic facts about plurisubharmonic functions, invariant metrics
and some standard theorems from complex analysis. The central tool in our proof
is the pluricomplex Green's function alluded to in the previous paragraph.

We will also give several applications that illustrate the scope of our
techniques.
A case in point is the situation of equality in \eqref{E:main}, which can be 
studied using our techniques in the case when 
$X = \disk$. This has been studied by
Nokrane and Ransford \cite[Theorem~1.2]{ransford2001schwarz} and our analysis is
reminiscent of theirs but simpler. 

\begin{theorem}\label{T:eqdsk}
  Let $f:\disk \to \sympoly$ be a holomorphic function such that 
  \begin{equation}
    \label{E:equal}
    \mathcal{H}^n_{\mathcal{M}_\disk}\big(\pi^{-1}(f(x)), \pi^{-1}(f(0)\big) =
    \mathcal{M}_\disk(x,0),
  \end{equation}    
  for $x \in U$ and $U \subset \disk$ a non-empty open subset. Then we can
  find an automorphism of $\disk$, say $g$, such that $\symf{g} \circ f$
  is the $n$-th root multi-function, i.e., the map
  \[
    z \mapsto \pi(\zeta_1(z),\dots,\zeta_n(z)),
  \] 
  where $\zeta_1(z),\dots,\zeta_n(z) \in \disk$ are the $n$-th roots of $z$.
\end{theorem}
\noindent Here $\sympoly$ is the symmetrized polydisk and $\pi$ is the
map whose coordinates are the elementary symmetric polynomials; see Section~
\ref{s:sym} for precise definitions including that of $\symf{g}$.

As another application, we shall also use our techniques to give a Schwarz lemma
for quasi-balanced domains that extends the well-known Schwarz lemma for balanced domains 
(Result~\ref{R:schwarz}); see Theorem~\ref{T:schwarz}. Using this lemma, we
shall then sketch a straightforward proof of a version of Schwarz lemma for the
spectral unit ball originally proved by Bharali \cite{bharali2007interpolation}.

\subsection*{Organization} Section~2 contains a brief treatment of all the
tools required in our proofs. We present our Schwarz lemma for quasi-balanced 
domains in Section~\ref{S:quasi}. The proofs of Theorems~\ref{T:main} and 
\ref{T:eqdsk} are contained in Section~\ref{S:proofs}. Finally, we briefly
sketch the proof of a version of the Schwarz lemma for the spectral unit ball 
in Section~\ref{S:spectral}.

\subsection*{Notation}
We will use $\disk$ to denote the unit disk in the
complex plane. The space of holomorphic mappings from a complex manifold $X$
into a complex manifold $Y$ will be denoted $\hol(X,Y)$. We use $|\bcdot|$ for the usual Euclidean norm
in $\C^n$, irrespective of the dimension. All manifolds will be assumed to be
connected. All other notations used will be introduced in Section~\ref{S:tools}.

\subsection*{Acknowledgments} I would like to thank Dr. Pranav Haridas for many
useful discussions. Prof. Aprameyan Parthasarathy and Dr. Vikramjeet Chandel
read parts of this manuscript and made useful suggestions that have greatly
improved the manuscript. Prof. G.P. Balakumar patiently answered  many
questions, both trivial and hard, that were useful in the proofs. I
would also like to thank Prof. Autumn Kent for bringing Stanton's paper 
\cite{stanton1975bounded} to my attention. I am grateful to the anonymous
referee for her/his careful reading of our paper and the many
useful suggestions that have improved our exposition.

\section{Tools}\label{S:tools}

\subsection{The pluricomplex Green's function}

In this section, we define and prove basic facts about an extremal function
defined using plurisubharmonic functions. Our treatment is from 
\cite[p.~184]{kobayashi1998} where the definition is attributed to Klimek 
\cite{klimek1985extremal}. The paper by Demailly 
\cite{demailly1987measures} contains further properties of this function.

\begin{definition}\label{D:extremal}
  Let $X$ be a complex manifold. Fix $z_0 \in X$ and define the \textbf{extremal
  function}
\begin{equation}\label{E:lambda}
  \lambda_X(z,z_0) := \sup\{\phi(z): \phi \in P_X(z_0)\},
\end{equation}
where $P_X(z_0)$ is the collection of functions $\phi$ on $X$ that satisfy:

\begin{enumerate}
  \item $\phi$ is upper semi-continuous,
  \item $0 \leq \phi < 1$,
  \item $\log \phi$ is plurisubharmonic on $X$,
  \item $\phi(z_0) = 0$,
  \item for any coordinate system $z = (z_1,\dots,z_n)$ with origin at $z_0$,
  the quantity $\frac{\phi(x)}{|z(x)|}$ is bounded above in a neighbourhood
  of $z_0$.
\end{enumerate}
\end{definition}

\begin{remark}\label{R:nonempty}
  In the above definition, functions that are identically $-\infty$ are
  considered to be plurisubharmonic whence the function that is identically $0$
  is an element of $P_X(z_0)$. So the collection $P_X(z_0)$ is always
  non-empty.
\end{remark}

\begin{remark}
  The function $\log \lambda_X(z,z_0)$ is known in the literature as the 
  \textbf{pluricomplex Green's function with a logarithmic pole at $z_0$}. The
  pluricomplex Green's function is well-studied and is at the heart of
  many deep results (see \cite{klimek1991pluri} and the papers cited in the
  introduction for a small sample). For our
  purposes, the function $\lambda_X$---which we will refer to throughout this
  paper as the \emph{extremal function}---is more convenient to work with. 
\end{remark}

\begin{remark}
  If $D \subset \C^n$ is a bounded domain then for each $z_0 \in D$, the function
  $|z - z_0| \in P_D(z_0)$. Therefore $\lambda_D(z,z_0)
  > 0 \ \ \forall z \in D \setminus \{ z_0 \}$.
\end{remark}

\begin{lemma}\label{L:dist}
  Let $X$ and $Y$  be complex manifolds and let $f:X \to Y$ be holomorphic. 
  Then
  \[
    \lambda_{Y}(f(x),f(z_0)) \leq \lambda_{X}(x,z_0).
  \]
\end{lemma}

\begin{proof}
  It suffices to show that if $\phi \in P_Y(f(z_0))$ then $\phi \circ f \in
  P_X(z_0)$. Only the final condition in the definition of $P_X(z_0)$
  needs to be checked. For a coordinate system $z = (z_1,\dots,z_n)$
  around $z_0$ and $w = (w_1,\dots,w_n)$ around $f(z_0)$, we have
  \begin{align*}
    \log \phi \circ f(x) - \log (|z(x)|) &= \log \phi \circ f(x)
    - \log |w(f(x))|\\ 
     &+ \log \frac{|w(f(x))|}{|z(x)|}.
  \end{align*}
  The expression on the right hand side is clearly bounded above in a
  neighbourhood of $z_0$ and
  we are done.
\end{proof}

We need a version of Schwarz lemma for subharmonic functions proved by Sibony in
order to compute the extremal function for the unit disk $\disk$.

\begin{lemma}[Sibony \cite{sibony1981hyperbolic}]\label{L:sibony}
  Let $u$ be an upper semi-continuous function on $\disk$ such that
  \begin{enumerate}
    \item $\log u$ is subharmonic,
    \item $\frac{u(z)}{|z|^2}$ is bounded on $\disk^*$,
    \item $0 \leq u < 1$ on $\disk$.
  \end{enumerate}
  Then $u(z) \leq |z|^2 \ \forall z \in \disk$. If $u(z_0) = |z_0|^2$ for some
  $z_0 \in \disk, z_0 \neq 0$, then $u(z) \equiv |z|^2$.
\end{lemma}

\begin{lemma}\label{L:disk}
  The extremal function $\lambda_\disk(z,0) = |z|$.
\end{lemma}

\begin{proof}
  Clearly $\lambda_D(z,0) \geq |z|$. Conversely, if $\phi \in P_\disk(z,0)$ then
  $\phi^2$ is subharmonic and $\frac{\phi^2(z)}{|z|^2}$ is bounded above on $\disk^*$ by the final condition in the definition of $P_\disk(z_0)$.
  This means that $\phi(z) \leq |z|$ by Lemma~\ref{L:sibony} and we are done.
\end{proof}

\subsection{The Möbius pseudodistance}

We now define the Möbius pseudodistance of a complex manifold $X$ and
prove some of its key properties. 

\begin{definition}
  Let $X$ be a complex manifold. We define the \textbf{Möbius pseudodistance
  on $X$} to be 
  \[
    \mathcal{M}_X(z_1, z_2) := \sup\{|f(z_1)|: f \in \hol(X,\disk), f(z_2) = 0
    \} \ \ \forall z_1, z_2 \in X.
  \]
\end{definition}

\begin{remark}
  Observe that the above definition is analogous to that of the Carathéodory
  pseudodistance except that we use the Möbius distance of $\disk$ in the
  definition instead of the Poincaré distance. The proof that the above
  definition actually gives a pseudodistance follows along the same lines as
  that for the Carathéodory pseudodistance. As expected, holomorphic mappings
  are distance decreasing under this pseudodistance and biholomorphisms are
  isometries. It is also clear that if $C_X$ denotes the Carathéodory
  pseudodistance on $X$ then $\tanh C_X = \mathcal{M}_X$. See \cite[Chapter~2]
  {jarnicki2013invariant} for details.
\end{remark}

\begin{remark}
  It follows from Lemma~\ref{L:disk} that
  \[
    \lambda_\disk(z,z_0) = \mathcal{M}_\disk(z,z_0).
  \]
\end{remark}
\begin{remark}\label{R:mobball}
  Let $B(a,r)$ be the ball of radius $r$ centred at point $a \in \C^n$. Then
  \[
    \mathcal{M}(z, a) = \frac{|z - a|}{r} \ \ \forall z \in B(a,r).
  \]
\end{remark}

\begin{remark}\label{R:mobpoly}
  If $D_i$ are disks in the plane then for $(z_1, \dots,z_n)
  \in D_1 \times \dots \times D_n$, we have
  \[
    \mathcal{M}_{D_1 \times \dots \times D_n}\big((z_1,\dots,z_n),(a_1,\dots,a_n)\big) = \max_i \mathcal{M}_{D_i}(z_i,a_i).
  \]
\end{remark}

\begin{remark}\label{R:cont}
  Using the Remark~\ref{R:mobball} and the fact $\mathcal{M}_X$ is distance decreasing
  under the inclusion map, one easily shows that $\mathcal{M}_X$ is continuous on
  $X\times X$; see \cite[Proposition~2.6.1]{jarnicki2013invariant}.

\end{remark}

\begin{definition}
  We say that the complex manifold $X$ is \textbf{Carathéodory hyberbolic} if
  $\mathcal{M}_X$ is a distance.
\end{definition}

\begin{remark}
   Bounded domains are Carathéodory hyperbolic. This follows from
   the observation that if $z, w \in D, z \neq w$, then some coordinate
   projection is a bounded holomorphic function that separates $z$ and $w$. 
\end{remark}

\begin{lemma}
  Let $X$ be a complex manifold. Then
  \[
    0 \leq \mathcal{M}_X < 1.
  \]
\end{lemma}

\begin{proof}
  Let $z,w \in X$ be such that $\mathcal{M}_X(z,w) = 1$. Then by the very
  definition of $\mathcal{M}_X$, we can find a sequence of holomorphic functions
  $f_n : X \to \disk$ such that $f_n(w) = 0$ and $|f_n(z)| \to 1$. By Montel's
  theorem, $\hol(X,\disk)$ is a normal family. This means that some subsequence
  of $f_n$ must converge in  the compact-open topology to a holomorphic map
  $f:X \to \disk$. But this is absurd as $|f_n(z)| \to 1$.
\end{proof}

The next theorem gives the crucial link between the Möbius
pseudodistance of a complex manifold and its extremal
function. This link
is the central tool used in the proof of Theorem~\ref{T:main}.

\begin{proposition}\label{P:sub}
  Let $X$ be a complex manifold. Then for a fixed $z_0
  \in X$, the function $\mathcal{M}_X(\bcdot,z_0)$ is plurisubharmonic. In fact,
  $\log\mathcal{M}_X(\bcdot, z_0)$ is plurisubharmonic. 
\end{proposition} 

\begin{proof}
  From Remark~\ref{R:cont}, $\mathcal{M}_X$ is a continuous function. 
  The fact that $\mathcal{M}_X(\bcdot,z_0)$ is plurisubharmonic is now
  straightforward from the fact that $\mathcal{M}_X(\bcdot,z_0)$ is continuous and a 
  supremum of plurisubharmonic functions. The same argument also shows that the
  function $\log \mathcal{M}_X(\bcdot, z_0)$ is plurisubharmonic. 
\end{proof}

\begin{remark}
  It is now straightforward to prove that for any $z_0 \in X$, the function $
  \mathcal{M}(\bcdot,z_0) \in P_X(z_0)$. Thus, $\mathcal{M}_X(\bcdot,z_0) \leq
  \lambda_X(\bcdot,z_0)$.
\end{remark}




\subsection{The Minkowski function of a quasi-balanced domain}

Let $p_1,p_2,\dots,p_n$ be relatively prime positive integers. We say
that a domain $D\subset \C^n$ is \textbf{$(p_1,p_2,\dots,p_n)$-balanced
  (quasi-balanced)} if
\begin{align*}
  \lambda \bullet z \in D \ \ \forall
  \lambda\in \overline{\D}, \forall z \in D,
\end{align*}
where for $z = (z_1, z_2, \ldots, z_n) \in D$, we define
$ \lambda \bullet z := (\lambda^{p_1}z_1,\lambda^{p_2}z_2,
\dots,\lambda^{p_n}z_n)$.  If $p_1 = p_2 = \dots = p_n = 1$ above,
then we say $D$ is a \textbf{balanced domain} 
(balanced domains are also known as \emph{complete circular domains} in the
literature).

Given a $(p_1,p_2,\dots,p_n)$-balanced domain $D\subset \C^n$, we
define the Minkowski function $\h_D : \C^n \rightarrow \C$ by
\begin{align*}
  \h_D(z) := \inf\{t > 0:  \frac{1}{t} \bullet z \in D \}.
\end{align*}
Clearly $D = \{z \in \C^n: \h_D(z) < 1\}$ and
$\h_D(\lambda\bullet z) = |\lambda|\h_D(z)$.  This function was first studied
by Nikolov \cite{nikolovsymmetrized} (see also \cite{bharali2006balanced}).
It turns out that
$\h_D$ is plurisubharmonic if and only if $D$ is additionally pseudoconvex; see 
\cite[Lemma~2.3]{bharali2006balanced}.

Section~2.2 of 
\cite{jarnicki2013invariant} contains an extensive treatment of the properties
of the Minkowski function of both balanced and quasi-balanced domains.

\subsection{The symmetric product of a Riemann surface}\label{s:sym}

Let $X$ be a Riemann surface. Given $(x_1,\dots, x_n) \in X^n$, we denote by
$\langle x_1, \dots, x_n \rangle$ the image in the quotient topological
space $\Sym^n(X) := X^n/S_n$ under the
$S_n$-action on $X^n$ that permutes the entries of $(x_1,\dots, x_n)$. We will
also abbreviate the element 
\[
  \langle \underbrace{z_1,\dots,z_1}_\text{$\mu_1$-times}, 
  \underbrace{z_2,\dots,z_2}_\text{$\mu_2$-times},\dots, \underbrace{z_k,\dots,z_k
  }_\text{$\mu_k$-times}\rangle, \mu_1 + \dots + \mu_k = n,
\]
by
\[
  \langle z_1;\mu_1, \dots,z_k;\mu_k \rangle. 
\]

The map
\[
  X^n\ni (x_1,\dots, x_n)\,\longmapsto\,\langle x_1,\dots, x_n\rangle
  \; \;\forall (x_1,\dots, x_n)\in X^n
\]
will be denoted by $\symq^n$. We shall drop the superscript when
there is no ambiguity. 
It is easy to see that there is a natural complex structure on
$\Sym^n(X)$ that makes it a complex manifold of dimension $n$ (see
below). With this
complex structure, the map $\symq$ is a branched proper holomorphic
mappings whose set of critical points is
\[
  \{(z_1,\dots,z_n) \in X^n : z_i = z_j \ \text{for some} \ i \neq j\}. 
\]

The symmetrized polydisk $\sympoly$ is a quasi-balanced domain in
$\C^n$ with weights $(1,2,\dots,n)$
defined using the elementary symmetric polynomials as follows. Let $\sigma_j$,
$j = 1,\dots, n$, denote the elementary symmetric polynomial of degree $j$ in
$n$ indeterminates. The map $\pi^{(n)} : \C^n\to \C^n$ is
defined as:
\begin{align*}
  \pi^{(n)}(z_1,\dots, z_n)\,:=\,\big(\sigma_1(z_1,\dots, z_n), \sigma_2(z_1,\dots, z_n),
  \dots, \sigma_n(z_1,\dots, z_n)\big),&
  \\ (z_1,\dots, z_n)\in \C^n.&
\end{align*}
Again, we shall drop the superscript when there is no scope for confusion.

The \textbf{symmetrized polydisk}, $\mathbb{G}^n$, is defined as
$\mathbb{G}^n := \pi(\disk^n)$. It is easy to see that $\mathbb{G}^n$ is
a $(1,2,\dots,n)$-balanced domain in $\C^n$, whence $\mathbb{G}^n$ is a
holomorphic embedding of the
$n$-fold symmetric product of $\disk$ into $\C^n$. It is also easy to see that the 
Minkowski functional of $\sympoly$ is given by
\[
  \h_{\sympoly}(z_1,\dots,z_n) := \max\{|\lambda_1|,\dots,|\lambda_n|: \pi^{(n)}
  (\lambda_1,\dots,\lambda_n) = (z_1,\dots,z_n)\}.
\]
The above formula implies that $\sympoly$ is pseudoconvex. This also follows
by appealing to the fact that the proper image of a pseudoconvex domain is
pseudoconvex. This automatically means that $\sympoly$ is a domain of
holomorphy.

We now give a brief description of the complex structure on the topological
space $\Sym^n(X)$ when $X$ is Riemann surface.
Given subsets $V_j\subseteq X$ that are open, let us write:
  \[
    \langle V_1,\dots, V_n \rangle := \left\{\langle x_1,\dots x_n\rangle : x_j \in V_j,
    \; j = 1, \dots, n\right\}
  \]
  The set $\langle V_1,\dots, V_n \rangle$ is an open subset of $X^n_
  \textsf{sym}$ by the defining property of the quotient topology. Given a point
  $p \in \Sym^n(X)$, $p = \langle p_1,\dots p_n \rangle$,
  choose a holomorphic chart $(U_j, \varphi_j)$ of $X$ at $p_j$,
  $j = 1,\dots, n$, such that
 \[
    U_j\cap U_k = \emptyset \ \ \text{if $p_j\neq p_k$} \qquad \text{and}
    \qquad U_j = U_k \ \ \text{if $p_j = p_k$}.
 \]
 The above choice of local charts ensures that the map  $\varPsi_p: \langle
 U_1,\dots, U_n\rangle\to \C^n$ given by
 \[
   \varPsi_p: \langle x_1,\dots, x_n\rangle\,\longmapsto\,
   \big((\varphi_1(x_1),\dots,
    \varphi_n(x_n)),\dots,(\varphi_1(x_1),\dots,\varphi_n(x_n))\big)
 \]
is a homeomorphism. This follows from the Fundamental Theorem of Algebra. The
collection of all such charts $(\langle U_1,\dots, U_n\rangle, \varPsi_p)$
produces a holomorphic atlas on $\Sym^n(X)$. The following lemma is easy to
prove and we omit the proof.

\begin{lemma}\label{L:subvariety}
  Let $X$ be a Riemann surface and for $1 \leq k < n$, define
  \[
    V_k := \{ \langle z_1, \dots, z_n \rangle \in \Sym^n(X): \text{the set $
    \{z_1,\dots,z_n\}$ has \emph{precisely} $k$ elements} \}
  \]
  Then $V_k$ is an analytic subvariety of $\Sym^n(X)$.
\end{lemma}

The book \cite{jarnicki2013invariant} contains an exhaustive account of the
various properties of the symmetrized polydisk. The book 
\cite{whitney1972complex} is the canonical reference for the symmetric product
in general.

\section{A Schwarz lemma for quasi-balanced domains}\label{S:quasi}

The following version of Schwarz lemma for balanced domains is well-known. This
version follows easily from the fact that holomorphic maps contract under the
Lempert function and the relationship between the Lempert function and the
Minkowski function of a balanced pseudoconvex domain.

\begin{result}[Proposition~3.1.1 of \cite{jarnicki2013invariant}]
\label{R:schwarz}
  Let $D_1 \subset \C^m$ and $D_2 \subset \C^n$ be balanced pseudoconvex domains
  with Minkowski functions $\h_1$ and $\h_2$, respectively. Then given any
  holomorphic map $f: D_1 \to D_2$ with $f(0) = 0$, we have 
  \[
    \h_2(f(z)) \leq \h_1(z).
  \]
\end{result}

We will now prove an analogue of the above result for quasi-balanced domains.

\begin{theorem}\label{T:schwarzb}
  Let $D$ be a $(p_1,\dots,p_n)$-balanced pseudoconvex domain with 
  highest weight $p_n$. Then 
  \[
    \h_D^{p_n}(z) \leq \lambda_D(z,0) \leq \h_D(z).
  \]
\end{theorem}

\begin{proof}
  First observe that the pseudoconvexity of $D$ ensures that $\log
  \h_D$ is plurisubharmonic. Fix $0 < \eps < 1$ and consider the set 
  \[
    K := \{w \in D: \h_D(w) = \eps\}.
  \]
  Note that $0 \not \in \overline{K}$. Observe that for any $z \in D$ 
  such that $0 < \h_D(z) < \eps$, we can find $0 < t < 1$ such that for some $z'
  \in K$, we have $t \bullet z' = z$. Hence $\h_D^{p_n}(z) = t^{p_n} \h_D^{p_n}(z')$. As
  $0 \not \in \overline{K}$, we can trivially write the inequality
  \[
    \h_D^{p_n}(w) \leq C |w| \ \forall w \in K,
  \]
  for some $C > 0$, suitably large. It is also easy to see that $|z| \geq t^
  {p_n} |z'|$. Thus, 

  \[
    \h_D^{p_n}(z) = t^{p_n}\h_D^{p_n}(z') \leq Ct^{p_n} |z'| \leq C
     |z|.
  \]
  Therefore $\frac{\h_D^{p_n}(z)}{|z|}$ is bounded in a
  neighbourhood of $0$. Thus, $\h_D^{p_n} \in P_D(z,0)$ whence
  $\h^{p_n}(z,0) \leq \lambda_D(z,0)$. This inequality is obviously also true if $\h_D(z_0) = 0$.
  
  Now for a fixed $z$ with $\h_D(z) \neq 0$, consider the map
  \[
    \phi: \disk \ni \lambda \mapsto \lambda \bullet \frac{z}{\h_D(z)} \in D.
  \]
  We then have 
  \[
    \lambda_D(z,0) \leq \lambda_\disk(\h_D(z),0) = \h_D(z).
  \]
  On the other hand, if $\h_D(z) = 0$, for each $n \in \Z_+$, the element $
  n \bullet z \in D$. We repeat the above argument with the element $n \bullet
  z$ instead of $\frac{z}{\h_D(z)}$. It is clear that $\lambda_D(z,0) \leq 1/n$.
  This proves that $\lambda_D(z,0) \leq \h_D(z)$ and we are done.
\end{proof}

The above theorem yields the following analogue of Schwarz lemma for
pseudoconvex quasi-balanced domains.
\begin{theorem}[Schwarz Lemma] \label{T:schwarz}
  Let $D_1 \subset \C^n$ and $D_2 \subset \C^m $ be pseudoconvex
  quasi-balanced domains. If $f:D_1 \to D_2$ is holomorphic and $f(0) =
  0$ then
  \[
    \h_{D_2}^p(f(z),0) \leq \h_D(z,0) \ \forall z \in D,
  \]
  where $p$ is highest weight of the quasi-balanced domain $D_2$
\end{theorem}

\begin{remark}
  The above theorem subsumes Result~\ref{R:schwarz}.
  See \cite[Theorem~1.6]{bharali2006balanced} for a proof of the above Schwarz
  lemma using the Lempert function instead of the extremal function. 
\end{remark}

The following is a version of Schwarz lemma that follows from Theorem~
\ref{T:schwarz}. This result was proved by 
Ransford--Nokrane \cite{ransford2001schwarz} in a formulation involving
algebroid multi-functions. 

\begin{theorem}
  Let $f:\disk \to \sympoly$ be holomorphic with $f(0) = 0$ and
  $f(z) = \pi(\lambda_1,\dots,\lambda_n)$. Then
  \[
    \max\{|\lambda_1|,\dots,|\lambda_n|\} \leq |z|^{1/n}.
  \]
\end{theorem}

\begin{remark}
  As alluded to in the introduction, the above observation was the impetus for
  this paper.
\end{remark}

\section{Proofs of the main results}\label{S:proofs}

Our strategy is to establish that the function $\mathcal{H}^n_{\mathcal{M}_X}
(\symq^{-1}(z), \symq^{-1}(z_0))$ (see \eqref{E:main}), is intimately related to
the extremal function of $\Sym^n(X)$ via a function $h_1$ which we will define
below.

Let $X$ be a Riemann surface  and fix $z_0 \in \Sym^n
(X)$. Define the function $h_1:\Sym^n(X) \to [0,1)$ by
\begin{equation}\label{E:h1}
  h_1(z) := \max\left(\max_i\prod_j\mathcal{M}_X(z_i, a_j), \max_i\prod_j
  \mathcal{M}_X(z_j, a_i) \right), 
\end{equation}
where $z_0 = \langle a_1,\dots,a_n \rangle$ is a fixed point and we have
written $z$ as $\langle z_1,\dots,z_n \rangle$. We also define the function
$h:\Sym^n(X)\to [0,1)$ by
\begin{equation}\label{E:h}
  h(z) := \mathcal{H}_{\mathcal{M}_X}(\symq^{-1}(z), \symq^{-1}(z_0)) = \mathcal{H}_{
  \mathcal{M}_X}(\{z_1,\dots,z_n\}, \{a_1,\dots,a_n\}) ,
\end{equation}
where $\symq^{-1}(z)$ and $\symq^{-1}(z_0)$ are defined as in Theorem~\ref{T:main}. Observe that from the very definitions, 
we have 
\begin{equation}\label{E:hh1}
  h^n(z) \leq h_1(z) \ \forall z \in \Sym^n(X).
\end{equation}

Our proof of the main theorem hinges on the following
theorem combined with Lemma~\ref{L:dist} and the fact that 
$\lambda_\disk(x,x_0) = \mathcal{M}_\disk(x,x_0)$. We, once again, emphasize 
that the Riemann surface $X$ is \emph{arbitrary} and in view of this,
Remark~\ref{R:nonempty} is pertinent in what follows.

\begin{theorem}\label{T:extremal}
   Let $V$ be the set of critical values of the map $\symq:X^n \to
   \Sym^n(X)$. For each $z_0 \in \Sym^n(X) \setminus V$, defining $h_1$ as in
   \eqref{E:h1}, we have
  \begin{equation} \label{E:klimineq}
    h_1 \in P_{\Sym^n(X)}(z_0). 
  \end{equation}
\end{theorem}

\begin{proof}
  From the very definition, $h_1$ is continuous, $h_1
  (z_0) = 0$ and $0 \leq h_1 < 1$. 
  We first show that the function $\log h_1$ is plurisubharmonic on
  $\Sym^n(X)$. Fix $z \in \Sym^n(X)\setminus V, z = \langle
  z_1,\dots,z_n \rangle$. Let $(U, \psi)$ be a coordinate chart around $z$ such
  that $\psi(z) = 0$. We can find an open set $B \subset U$ such that: 
  \begin{enumerate}
    \item The map $\psi|_B$ is a biholomorphism onto a ball $B(0,r)$,
    \item We can find  an inverse $(\widetilde{\pi}_1,\dots,\widetilde{\pi}_n)$
    of $\symq$ defined on $B$ such that $\widetilde{\pi}_i(z) = z_i$.
  \end{enumerate}

  For $y \in B$, we can write
  \[
    h_1(y) = \max\left(\max_i\prod_j\mathcal{M}_X(\widetilde{\pi}_i(y), a_j),
    \max_i\prod_j\mathcal{M}_X(\widetilde{\pi}_j(y), a_i) \right).
  \]
  Now Proposition~\ref{P:sub}, together with basic properties of plurisubharmonic
  functions, shows that $\log h_1$ is plurisubharmonic on $\Sym^n(X)
  \setminus V$. By Riemann's removable singularities theorem for plurisubharmonic  functions (\cite
  [Theorem~3, p.~178]{gunning1990holomorphic}), the function $h_1$ extends to be
  a plurisubharmonic function on $\Sym^n(X)$.

  It remains to show that the final condition in the definition of $P_X(z_0)$
  is satisfied by $h_1$. Let $(U, \psi)$ be any coordinate chart
  around $z_0$ such that $\psi(z_0) = 0$. Choose $B$ and $
  \widetilde{\pi}_i$ as before. Let $D_i
  \subset X$ be open pairwise disjoint coordinate disks that contain $a_i$. By
  continuity, shrinking $B$ if necessary, we can assume $(
  \widetilde{\pi}_1,\dots, \widetilde{\pi}_n)(B) \subset
  D_1 \times \dots \times D_n$. By the distance decreasing
  property of the Möbius pseudodistance and Remark~\ref{R:mobpoly}, we now have
  \[
    \max_i\mathcal{M}_{D_i}(\widetilde{\pi}_i(z), a_i) \leq \mathcal{M}_B
    (z,z_0) \ \ \forall z \in B.
  \]
  From Remark~\ref{R:mobball}, $\mathcal{M}_B(z,z_0) = \frac{|\psi(z)|}{r}$. The
  above equation, combined with the fact that $\mathcal{M}_X \leq \mathcal{M}_
  {D_i}$, therefore shows
  \[
    \max_i \mathcal{M}_X(\widetilde{\pi}_i(z), a_i) \leq \frac{|\psi(z)|}{r} \
    \ \forall z \in B.
  \]
  From the very definition of $h_1$, it is now follows that
  \[
    h_1(z) \leq \frac{|\psi(z)|}{r} \ \ \forall z \in B.
  \]
  The function $h_1$ satisfies all the conditions required for it to be an
  element of $P_X(z_0)$ and we are done.

\end{proof}

\begin{remark}
  It is not hard to see that $\lambda_X(\bcdot,z_0) \in P_X(z_0)$ (see 
  \cite[Corollary~1.3]{klimek1985extremal}). Therefore, in the definition of 
  $h_1$, we might as well have used the function $\lambda_X$ instead of the
  function $\mathcal{M}_X$ and the same proof \emph{mutatis mutandis} would
   show that the modified function is in $P_X(z_0)$ as well.
\end{remark}

The following corollary is immediate from \eqref{E:hh1} and Theorem~
\ref{T:extremal}.
\begin{corollary}\label{C:haus}
  For each $z_0 \in \Sym^n(X) \setminus V$
  \[
    h^n(z) \leq h_1(z) \leq \lambda_{\Sym^n(X)}(z,z_0).
  \]
\end{corollary}

Before we come to the proof of Theorem~\ref{T:main}, we need one final lemma.

\begin{lemma}\label{L:lift}
  With the same notation as Theorem~\ref{T:main}, let $1 \leq k \leq n$ be the
  highest integer such that for some $x^0 \in \disk$, writing $f(x^0) = \langle
  x^0_1,\dots,x^0_n \rangle$,  the set $\{x^0_1,\dots,x^0_n\}$ has $k$ elements.
  Then: 
  \begin{enumerate}
    \item Except for $x$ in a discrete set
          $E \subset \disk$, $f(x) = \langle x_1,\dots,x_n \rangle$ also has
          the  property that $\{x_1,\dots,x_n\}$ has $k$ elements;

    \item  For each $x \in \disk \setminus E$, we can find a disk $
           V_x \subset \disk \setminus E$ centred at $x$, holomorphic maps $
           \widetilde{f}_{x,1},\dots, \widetilde{f}_{x,k} : V_x \to X$ and
           positive integers $\mu_1, \dots, \mu_k$ whose sum is $n$ such that 
           \[
             f(y) = \left \langle 
             \widetilde{f}_{x,1}(y);\mu_1,\widetilde{f}_{x,2}(y);\mu_2,\dots,
             \widetilde{f}_{x,k}(y);\mu_k \right \rangle, \forall y \in V_x.
           \]
  \end{enumerate}
 
\end{lemma}

\begin{proof}
  Let $E \subset \disk$ be the set of all elements each $x \in \disk$ with the
  property that $f(x) = \langle x_1,\dots,x_n \rangle$ is such that $
  \{x_1,\dots,x_n\}$ has fewer than $k$ elements. By Lemma~\ref{L:subvariety},
  the collection of
  all points $w$ in $\Sym^n(X)$ with the property that, writing $w$ as
  $\langle w_1,\dots, w_n \rangle$, the set $\{w_1,\dots,w_n\}$ has
  fewer than $k$ elements is an  analytic subvariety of $\Sym^n(X)$ 
  (Lemma~\ref{L:subvariety}). If $E$ is an indiscrete set, it follows from the 
  principle of analytic continuation that $E = \disk$, a  contradiction.

  Now let $x \in \disk \setminus E$ and $f(x) = \langle x_1;\mu_1,\dots x_k;
  \mu_k \rangle, \mu_1 + \dots + \mu_k =n$. Let $U_i \subset X$ be pairwise
  disjoint coordinate disks
  centred at $x_i$. Then by continuity, we can find a disk $V_x \subset \disk
  \setminus E$ centred at $x$ such that 
  \[
    f(V_x) \subset \left\langle \underbrace{U_1, \dots, U_1}_{\mu_1-
    \text{times}},\dots, \underbrace{U_k, \dots, U_k}_{\mu_k-\text{times}}
    \right\rangle.
  \]
  As the $U_i$ are pairwise disjoint and for each $y \in V_x$ and writing $f(y) =
  \langle y_1,\dots,y_n \rangle$, the cardinality of $\{ y_1,\dots,y_n\}$
  is $k$, it is clear that we can define continuous
  maps $\widetilde{f}_{x,1}, \dots, \widetilde{f}_{x,k} : V \to X$ such that
   \[
    f(y) = \left \langle 
    \widetilde{f}_{x,1}(y);\mu_1,\widetilde{f}_{x,2}(y);\mu_2,\dots,
    \widetilde{f}_{x,k}(y);\mu_k \right \rangle, \forall y \in V_x.
   \]
  The fact that the maps $\widetilde{f}_{x,1}, \dots, \widetilde{f}_{x,k}$
  are holomorphic is a simple consequence of the way the complex structure on
  $\Sym^n(X)$ is defined.

\end{proof}

\subsection*{Proof of Theorem~\ref{T:main}} 
Let $E$ and $k$ be as in Lemma~\ref{L:lift}. For $x \in
\disk \setminus E$, we can find a disk $V_x$ and holomorphic maps $
\widetilde{f}_{x,1}, \dots, \widetilde{f}_{x,k}$ defined on $V_x$ satisfying
the conclusion of Lemma~\ref{L:lift}. We now define
$\hat{f}$ on $V_x$ by
\[
  V_x \ni y \mapsto \symq^k\left(\widetilde{f}_{x,1}(y), \dots, 
  \widetilde{f}_{x,k}(y)\right). 
\]
The above definition yields a holomorphic map $\hat{f}: \disk \setminus E
\to \Sym^k(X)$.
By Riemann's removable singularities theorem for subharmonic functions, 
$\lambda_{\disk \setminus E} \equiv \lambda_\disk|_{\disk
\setminus E}$. Fix $y \in \disk \setminus E$ and define the functions $
\hat{h}$
and $\hat{h}_1$ on $\Sym^k(X)$ with respect to the point $\hat{f}(y)$
and analogous to $h$ and $h_1$ (see \eqref{E:h1} and \eqref{E:h}), 
respectively. It follows from Corollary~
\ref{C:haus} and Lemma~\ref{L:dist} that
\[
  \hat{h}^n(\hat{f}(x)) \leq \hat{h}_1(\hat{f}(x)) \leq \mathcal{M}_\disk(x,
  y) \ \forall x \in \disk \setminus E.
\]
It is obvious that 
\[
  h_1(f(x)) \leq \hat{h}_1(\hat{f}(x)),
\] 
where $h_1$ is defined on $\Sym^n(X)$ with respect to the point $f(y)$.
Thus
\[
  h_1(x) \leq \mathcal{M}_\disk(x,y) \ \forall x \in \disk\setminus E.
\]
From \eqref{E:h1}, it is clear that if we view $h_1$ as a 
function of both $x$ and $y$, it is continuous on $\disk \times \disk$. So is
the function $\mathcal{M}_\disk(x,y)$. This combined with \eqref{E:hh1}
delivers the theorem.

\subsection*{Proof of Theorem~\ref{T:eqdsk}}
  In this proof, we shall tacitly identify $\Sym^n(\disk)$ with $\sympoly$
  without explicit mention. 
  Let $z_0 := f(0) = \langle a_1, \dots, a_n \rangle$ and consider the
  functions $h$ as before defined on $\Sym^n(\disk)$ with respect to the
  point $f(0)$. It is harmless to assume that $0 \not \in U$. 

  \noindent \textbf{Claim:} \emph{We can find an open disk $G \subset U$
  and a holomorphic function $F:G\to \disk$ such that for some $1 \leq j_0
  \leq n$, we have
  \[
    h^n(x) = \mathcal{M}_\disk(F(x), a_{j_0}) \ \ \forall x \in G.
  \]}

  \noindent \emph{Proof of claim:} 
  We adopt the same notation as Lemma~\ref{L:lift}. Choose $x^0 \in U \setminus E$.
  We have $h(x^0) = \mathcal{M}_\disk(\widetilde{f}_{x^0, i_0}(x^0),
  a_{j_0})$ for some choice of $1 \leq i_0, j_0 \leq k$ (the choice might not
  be unique). Let $i_0,\dots,i_l$ and $j_0,\dots,j_l$ be all the indices such
  that $h(x^0) = \mathcal{M}_\disk(\widetilde{f}_{x^0, i_m}(x^0), a_{j_m})$
  where $0 \leq m \leq l$.  We can find a disk $G \subset V_{x^0} \cap U$
  centered at $x^0$
  such that for each $x \in G$, $h(f(x))$ is one of the functions $\mathcal{M}_\disk(
  \widetilde{f}_{x^0, i_m}
  (x), a_{j_m}), 1 \leq m \leq l$. Define the sets
  \[
    E_m := \left\{x \in G: h(f(x)) = \mathcal{M}_\disk\left(\widetilde{f}_
    {x^0,i_m}(x), a_{j_m}\right)\right\},  0 \leq m \leq l.
  \]
  Each $E_m$ is a closed subset of $G$ and $\bigcup_{m=0}^l E_m = G$.
  Consequently, one of the sets $E_m$ has non-empty interior and we can 
  rename $G$ to be any disk contained in this $E_m$ and choose $F$ to be the
  corresponding $\widetilde{f}_{i_m}$.

  With the claim in hand, the proof of the theorem is not hard. Let $a_{j_0}$,
  $G$ and $F$ be as in the claim. We may assume that $a_{j_0} \not\in G$. Let
  $\phi \in \textsf{Aut}(\disk)$ be the automorphism that
  interchanges $0$ and $a_{j_0}$. We have $h^n(\symf{\phi} \circ f(x)) =
  \mathcal{M}_\disk(\phi \circ F(x),0))^n = |\phi \circ F(x)|^n \ \forall x \in
  G$. By hypothesis, this means that
  \[
    |\phi \circ F(x)|^n = |x| \ \forall x \in G.
  \]
  But any branch of $\sqrt[n]{\bcdot}$ on $G$ satisfies the
  above equation
  as well proving that for some $\theta$, $e^{i\theta}(\phi \circ F)$ is just
  some branch of the $\sqrt[n]{\bcdot}$. Let $\symf{\Theta}$ be the automorphism of $\sympoly$ associated to rotation by $e^{i\theta}$. 
  Replacing $f$ with $\symf{\Theta} \circ \symf{\phi} \circ f$, we
  may assume that $f|_G$ lifts over $\pi$ to a map into $\D^n$, one of whose
  components is a branch of $\sqrt[n]{\bcdot}$.

  Writing $f$ as $(f_1,\dots,f_n)$, consider the polynomial over $\hol(\disk)$
  \[
    P(x,y) := y^n + f_1(x) y^{n-1} + \dots + f_{n-1}(x)y^{n-1} + f_n(x)
  \] 
  From the conclusion of the preceding paragraph,  we can find a $n$-th root of
  unity $\zeta$  such that $P(x^n, x \zeta) \equiv 0$ on $G$. Consequently, $P
  (x^n, x \zeta) \equiv 0$ on $\disk$ by the identity theorem. If $\eta$ is any
  other $n$-th root of unity, we see that $P(x^n, x \eta \zeta) \equiv 0$ on
  $\disk$. Therefore $f(x) = \pi(\zeta_1(x),\dots,\zeta_n(x))$
  where $\sqrt[n]{x} = \{\zeta_1(x), \dots, \zeta_n(x)\}$. The theorem
  is proved with $g := e^{i\theta} \phi$.

\section{A Schwarz lemma for the spectral unit ball}\label{S:spectral}

In this section, we sketch a proof of a Schwarz lemma for the spectral
unit ball. This theorem was formulated and proved by Bharali 
\cite{bharali2007interpolation}. But as the ideas fit well with the main themes
of this article, we felt it is worthwhile to sketch a slightly different proof
here.

For $n \in \Z_+$, the spectral unit ball $\OM_n \subset \C^{n^2}$ is the
collection of all matrices $A \in M_n(\C)$ ($n \times n$ complex matrices) whose
spectrum $\sigma(A)$ is contained in $\disk$.  We have the following

\begin{proposition}
  The set $\OM_n$ is an unbounded balanced pseudoconvex domain with Minkowski 
  function given by the spectral radius $\rho$.
\end{proposition}

\begin{proof}
  That $\OM_n$ is balanced and that the spectral radius is the Minkowski
  function is easy to see from the
  definitions. We can define the holomorphic map $\Psi_n:M_n(\C) \to \C^n$ given
  by $M \mapsto \pi(\sigma(M))$. Observe that $\Psi_n^{-1}(\sympoly) = \OM_n$
  which shows that $\OM_n$ is a domain of holomorphy  (from that fact that
  $\sympoly$ is a domain of holomorphy and \cite[Theorem~2.5.14]{hormander}).
  Pseudoconvexity of $\OM_n$ now
  follows from the characterization of domains of holomorphy (see \cite
  [Section~2.6]{hormander}).
\end{proof}

\begin{remark}
  The above proposition shows that $\rho|_{\OM_n}$ is
  plurisubharmonic (see \cite[Appendix B.7.6]{jarnicki2013invariant}). This fact
  is usually proved in the literature using a theorem of 
  Vesentini \cite{vesentini1968subharmonic}.  
\end{remark}

\begin{definition}
  Given $A\in M_n(\C)$, we can write its minimal polynomial $\mathbf{M}_A$
  as
  \[
    \mathbf{M}_A(t)=\sum_{\lam\in\sigma(A)}(t-\lam)^{m(\lam)}.
  \]
  The \textbf{minimal Blaschke product corresponding to $A$} is defined by
  \begin{equation}\label{E:minimal}
    \mathbf{B}_A(t):=\prod_{\lambda\in\sigma(A)\subset\D}{\left(\frac{t-\lambda}
    {1 - \overline{\lambda}t}\right)^{m(\lambda)}}.
  \end{equation}
\end{definition}

Using the minimal Blaschke product corresponding to $A$, we can construct a
holomorphic map $\widetilde{A}: \OM_n \to \OM_n$ that takes $A$ to $0$. We
define
\[
  \widetilde{A}:B \mapsto \prod_{\lambda \in \sigma(A)}(\mathbb{I} - 
  \overline{\lambda}B)^{-m(\lambda)}(B - \lambda \mathbb{I})^{m(\lambda)},
\]
where $m(\lambda)$ is the multiplicity of the eigenvalue $\lambda$ in the
minimal polynomial of $A$. It can be shown that if $\sigma(B) = \{\lambda_1,
\dots, \lambda_n\}$ then $\sigma(\widetilde{A}(B)) = \left\{\mathbf{B}_A
(\lambda_1), \dots, \mathbf{B}_A(\lambda_n)\right\}$.
If $F:\disk \to \OM_n$ is holomorphic such that $F(z) = A$ and $F(w) = B$ then
$\widetilde{A} \circ F$ takes $A$ to $0$ and $\widetilde{B} \circ F$ takes
$B$ to $0$.
The following result is immediate from the Schwarz lemma for
balanced domains (Result~\ref{R:schwarz}). 

\begin{result}[Bharali, Theorem~1.5 of \cite{bharali2007interpolation}] 
  Let $f:\disk \to \OM_n$ be holomorphic. Then for $z, w \in \disk$, we have
  \begin{align*}
     \max & \left\{\max_{\lambda \in \sigma(f(w))}\prod_{\mu \in \sigma(f(z))}
     \mathcal{M}_\disk(\mu,\lambda)^{m(\mu)},\max_{\mu \in \sigma(f(z))}\prod_
     {\lambda
     \in
     \sigma(f(w))}\mathcal{M}_\disk(\mu,\lambda)^{m(\lambda)}\right\} \\
     & \hspace{3.5in}\leq \mathcal{M}_\disk(z,w),
  \end{align*}
  where $m(\mu)$ and $m(\lambda)$ denote the multiplicity of the eigenvalues
  $\mu$ and $\lambda$ in $\mathbf{M}_{f(z)}$ and $\mathbf{M}_{f(B)}$,
  respectively.
  
\end{result}

\bibliographystyle{amsalpha} 
\bibliography{schwarz}
\end{document}